\documentclass{gtart_h}

\def\ifplaintex{\expandafter\ifx\csname documentclass\endcsname\relax}

\ifplaintex 
\else
\fi

\def\gt{{\mathsurround=0pt\it $\cal G\mskip-2mu$eometry \&\ 
$\cal T\!\!$opology}}        

\def\gtm{{\mathsurround=0pt\it $\cal G\mskip-2mu$eometry \&\ 
$\cal T\!\!$opology $\cal M\mskip-1mu$onographs}}    

\def\gtp{{\mathsurround=0pt\it $\cal G\mskip-2mu$eometry \&\ 
$\cal T\!\!$opology $\cal P\!$ublications}}  

\def\recd{{\small Received:\qua\receiveddate\ifx\reviseddate\relax
\else\qquad Revised:\qua\reviseddate\fi\par}} 

\def\volumenumber#1{\def\thevolumenumber{#1}}
\def\volumeyear#1{\def\thevolumeyear{#1}}
\def\volumename#1{\def\thevolumename{#1}}
\def\papernumber#1{\def\thepapernumber{#1}}
\def\pagenumbers#1#2{\def\startpage{#1}\def\finishpage{#2}}
\def\published#1{\def\publishdate{#1}}
\def\received#1{\def\receiveddate{#1}}
\def\revised#1{\def\reviseddate{#1}}
\def\accepted#1{\def\accepteddate{#1}}

\def\coverauthors#1{\def\thecoverauthors{#1}}
\def\asciiauthors#1{\def\theasciiauthors{#1}}
\def\asciiaddress#1{\def\theasciiaddress{#1}}
\def\asciiemail#1{\def\theasciiemail{#1}}

\def\coverauthors#1{\def\thecoverauthors{#1}}

\let\\\par
\let\thevolumenumber\relax\let\thepapernumber\relax
\let\thevolumeyear\relax\let\startpage\relax
\let\finishpage\relax\let\publishdate\relax\let\receiveddate\relax
\let\reviseddate\relax\let\accepteddate\relax\let\theasciititle\relax
\let\theasciiauthors\relax\let\theasciiaddress\relax
\let\theasciiabstract\relax
\let\thecoverauthors\relax
\let\thecoverauthors\relax\let\theerratum\relax\let\theasciiemail\relax
\let\theshortauthors\relax\let\theshorttitle\relax

\def\startpage{1}\def\finishpage{15}\def\thepapernumber{77}

\volumenumber{2}
\volumename{Proceedings of the Kirbyfest}
\volumeyear{1999}

\long\def\maketitlep{   

\count0=\startpage

\gtm\nl        
{\small Volume \thevolumenumber: \thevolumename\nl 
\ifx\theerratum\relax\else Erratum \erratumnumber\nl\fi
Pages \startpage--\finishpage\nl}

\vglue 0.1truein   

{\parskip=0pt\leftskip 0pt plus 1fil\def\\{\par\smallskip}{\ifplaintex\large
\else\Large\fi\bf\thetitle}\par\medskip}   
\vglue 0.05truein 

{\parskip=0pt\leftskip 0pt plus 1fil\def\\{\par}{\sc\theauthors}
\par\medskip}%
 
\vglue 0.03truein 

{\small\leftskip 25pt\rightskip 25pt{\bf Abstract}\stdspace\theabstract

{\bf AMS Classification}\stdspace\theprimaryclass
\ifx\thesecondaryclass\relax\else; \thesecondaryclass\fi\par
{\bf Keywords}\stdspace \thekeywords\par}\vglue 7pt

}   

\font\phead=cmsl9 scaled 950
\font=cmsl9 scaled 1050
\font\pnum=cmbx10 scaled 913
\font\lnum=cmbx10 
\font\pfoot=cmsl9 scaled 950
\font=cmsl9 scaled 1050
\ifplaintex
\headline{\vbox to 0pt{\vskip -4.5mm\line{\small\phead\ifnum
\count0=\startpage ISSN 1464-8997 (on line)
1464-8989 (printed) \hfill {\pnum\folio}\else\ifodd\count0\def\\{ }%
\ifx\theshorttitle\relax\thetitle\else\theshorttitle\fi\hfill{\pnum\folio}
\else\def\\{ and }{\pnum\folio}\hfill\ifx\theshortauthors\relax\theauthors
\else\theshortauthors\fi\fi\fi}\vss}}
\footline{\vbox to 0pt{\vglue 0mm\line{\small\pfoot\ifnum\count0=\startpage
Published \publishdate:\qua\copyright\ \gtp\hfill\else
\gtm, Volume \thevolumenumber\ (\thevolumeyear)\hfill\fi}\vss
}}
\else
\makeatletter
\def\@oddhead{{\small\lhead\ifnum\count0=\startpage ISSN 1464-8997 (on line)
1464-8989 (printed) \hfill {\lnum\number\count0}\else\ifodd\count0
\def\\{ }\ifx\theshorttitle\relax \thetitle \else\theshorttitle\fi\hfill
{\lnum\number\count0}\else\def\\{ and }{\lnum\number\count0}
\hfill\ifx\theshortauthors\relax 
\theauthors\else\theshortauthors\fi\fi\fi}}\def\@evenhead{@oddhead}
\def\@oddfoot{\small\lfoot\ifnum\count0=\startpage Published \publishdate:\qua\copyright\ \gtp\hfill\else
\gtm, Volume \thevolumenumber\ (\thevolumeyear)\hfill\fi}
\def\@evenfoot{@oddfoot}
\makeatother
\fi

\let\maketitlepage\maketitlep
\let\makeshorttitle\maketitlepage
\let\maketitle\maketitlepage

\newwrite\gtoutfile
\long\gdef\makeheadfile{  
{\def\\{, }\def\s{ }
\immediate\openout\gtoutfile head.xxx
\immediate\write\gtoutfile{Proxy-for: \ifx\theasciiauthors\relax
\theauthors\else\theasciiauthors\fi\s<\ifx\theasciiemail\relax\theemail\else\theasciiemail\fi>}
\immediate\write\gtoutfile{\noexpand\\}
\immediate\write\gtoutfile{Authors: \ifx\theasciiauthors\relax
\theauthors\else\theasciiauthors\fi}
{\def\\{ }\immediate\write\gtoutfile{Title: \ifx\theasciititle\relax
\thetitle\else\theasciititle\fi}}
\immediate\write\gtoutfile{Subj-class: GT or SG, GR etc}
\immediate\write\gtoutfile{MSC-class: \theprimaryclass\ifx\thesecondaryclass\relax\else, \thesecondaryclass\fi}
\immediate\write\gtoutfile{Journal-ref: Geom. Topol. Monogr. \thevolumenumber\s
(\thevolumeyear) \startpage-\finishpage}
\immediate\write\gtoutfile{Comments: Published by Geometry and Topology Monographs at}
\immediate\write\gtoutfile{\s\s\s  http://www.maths.warwick.ac.uk/gt/GTMon\thevolumenumber/paper\thepapernumber.abs.html}
\immediate\write\gtoutfile{\noexpand\\}
\immediate\write\gtoutfile{}
\ifx\theasciiabstract\relax
\immediate\write\gtoutfile{\theabstract}\else
\immediate\write\gtoutfile{\theasciiabstract}\fi
\immediate\write\gtoutfile{}
\immediate\write\gtoutfile{\noexpand\\}
\immediate\write\gtoutfile{}
\immediate\closeout\gtoutfile}}  

\def\maketitlepage{\maketitlep\makeheadfile}
\let\makeshorttitle\maketitlepage
\let\maketitle\maketitlepage

\volumenumber{7}
\volumename{Proceedings of the Casson Fest}
\volumeyear{2004}
\papernumber{8}
\pagenumbers{205}{212}
\received{19 November 2003}
\revised{21 April 2004}
\accepted{20 April 2004}
\published{19 September 2004}

\usepackage{amsmath, amssymb}

\newtheorem{theorem}{Theorem}[section]
\newtheorem{lemma}[theorem]{Lemma}

\theoremstyle{definition}

\def\GL{\mbox{\rm{GL}}}
\def\PSL{\mbox{\rm{PSL}}}
\newcommand{\area}{\operatorname{area}}
\newcommand{\sign}{\operatorname{sign}}
\newcommand{\stab}{\operatorname{stab}}

\begin{document}

\title{Generalized Dedekind sums}
\authors{D\,D Long\\A\,W Reid}
\coverauthors{D\noexpand\thinspace D Long\\A\noexpand\thinspace W Reid}
\asciiauthors{DD Long\\AW Reid}
\address{Department of Mathematics, University of California\\Santa
Barbara, CA 93106, USA\\{\rm and}\\
Department of Mathematics, University of Texas\\Austin, TX 78712, USA}
\asciiaddress{Department of Mathematics, University of California\\Santa
Barbara, CA 93106, USA\\and\\Department of Mathematics, University of Texas\\Austin, TX 78712, USA}
\asciiemail{long@math.ucsb.edu, areid@math.utexas.edu}
\gtemail{\mailto{long@math.ucsb.edu}\qua{\rm and}\qua 
\mailto{areid@math.utexas.edu}}

\begin{abstract}
Classical Dedekind sums are connected to the modular group
through the construction of a (Dedekind) symbol on the cusp set of the
modular group.  In this paper we study generalizations of Dedekind
symbols and sums that can be associated to certain Fuchsian groups
uniformizing 1-punctured tori.
\end{abstract}

\keywords{Cusp set, Dedekind sum}

\primaryclass{11F20}
\secondaryclass{20H10}
\maketitle

\cl{\small\em Dedicated to Andrew Casson on the occasion of his 60th birthday}

\section{Introduction}

A classical and important construction which arises in many contexts 
is that of the {\em Dedekind sum} which is defined for coprime 
integers $a$ and $c$ by
$$s(a,b) = \sum_{k=1}^{|b|-1} \left(\!\!\left(\frac{k}{b}\right)\!\!\right)
  \left(\!\!\left(\frac{ka}{b}\right)\!\!\right)$$
where $((x)) = x - [x] -1/2$. Dedekind sums arise naturally in 
various topological settings, one of the most famous being 
Hirzebuch's description of $4s(b,a)$ as the signature defect of the 
Lens space $L(a,b)$ coming from Rademacher's cotangent formula
$$ s(a,b) = \frac{1}{4|b|}\sum_{k=1}^{|b|-1} 
  \cot\left(\frac{k\pi}{b}\right)
  \left(\!\!\left(\frac{ka\pi}{b}\right)\!\!\right),$$
as well as in Walker's formula for the generalized Casson invariant.

From the point of view of this note, it is the beautiful construction
in \cite{KM} of Dedekind sums based upon the classical modular group
$\PSL(2,\mathbb{Z})$ that is of interest.  We describe some of this briefly, 
as it is useful in the development of what follows. It is shown in
\cite{KM} that there exists a 2--cocycle 
$\epsilon\co\PSL(2, \mathbb{Z})\times \PSL(2,\mathbb{Z})\rightarrow \mathbb{Z}$
and a function $\phi\co\PSL(2,\mathbb{Z})\rightarrow \mathbb{Z}$
(the Rademacher $\phi$--function) which satisfy $\delta \phi = 3\epsilon$
(where $\delta$ is the coboundary operator). Furthermore, it
is shown in \cite{KM} that the function $\phi$ is closely related to the
Dedekind sums mentioned above. Namely, in \cite{KM} the authors define
a Dedekind symbol $S$ on $\mathbb{Q}\cup \infty$ which maps $\infty$ to $\infty$
and otherwise, $S(\frac{a}{c}) = 
\phi(M) + \chi(M)$ where $M\in \PSL(2,\mathbb{Z})$ satisfies 
$M(\infty) = \frac{a}{c}$ and $\chi$ is a function depending on the entries
of $M$ (see section 2.2). As pointed out in \cite[section 0.8]{KM}, the relationship
between $S$ and the Dedekind sum $s$ above is
$S(\frac{a}{c}) = 12 \sign(c)s(a,c)$.

For us, since $\mathbb{Q}\cup \infty$
coincides with the cusp set (that is the set of all parabolic fixed points)
of $\PSL(2,\mathbb{Z})$,
$S$ can be viewed as a function defined on the cusp set of $\PSL(2,\mathbb{Z})$.
In \cite{LR} it was shown that there exist finite coarea Fuchsian groups not
commensurable with the modular group but whose cusp set is precisely
$\mathbb{Q} \cup \infty$.  The purpose of this note is to show that these
groups give rise to very natural generalizations of Dedekind sums.

We begin by recalling briefly the construction of \cite{LR}. The
starting point of that paper was to take the two generator group
$\Delta(u^{2},2t)$ generated by elements $g_{1}$ and $g_{2}$
as below
$$g_1 = \left( \begin{array}{ccc} (-1+t)/\sqrt{-1+t-u^2}& &
  u^2/\sqrt{-1+t-u^2} \\
& & \\ 1/\sqrt{-1+t-u^2}& & 
  1/\sqrt{-1+t-u^2} \end{array} \right)$$
and
$$g_2 = \left( \begin{array}{ccc} u/\sqrt{-1+t-u^2}& &
  u/\sqrt{-1+t-u^2} \\
& & \\ 1/(u\sqrt{-1+t-u^2})& & 
  (t-u^2)/u\sqrt{-1+t-u^2} 
\end{array} \right)$$
where the parameters $u^{2}$ and $t$ are real and satisfy $t>u^{2}+1$.

One sees easily that in the hyperbolic plane, $g_1$ maps the directed
edge $\{-1,0\}$ to the directed edge $\{\infty, u^2\}$ and $g_2$
mapping $\{\infty,-1\}$ to $\{u^2, 0\}$, and moreover the commutator
$$g_1 g_2^{-1} g_1^{-1} g_2 =
  \left( \begin{array}{cc} -1 & -2t \\ 0& -1 \end{array} \right)$$ is
parabolic and generates the stabiliser of
infinity. It follows that ${\bf H}^{2}/\Delta(u^{2},2t)$ is a
complete finite-area once-punctured torus. This family includes a
modular torus as $\Delta(1,6)$, as well as other arithmetic once-punctured
tori, and if $u^{2}$ and $t$ are chosen to be rational the
set of cusps of these groups must be a subset of $\mathbb{Q}\cup \infty$. In
the arithmetic cases, the cusp set is precisely $\mathbb{Q}\cup \infty$,
although this is not always the case for rational pairs
$(u^{2},2t)$. (See \cite{LR}).

Despite the apparently complicated nature of the entries in these
matrices because of the presence of square roots, an easy computation
shows that if one considers $G = \ker \{ \Delta \rightarrow \mathbb{Z}/2
\oplus \mathbb{Z}/2\}$, then the trace-field of $G$, and hence the
invariant trace-field of $\Delta(u^2,2t)$
is the field $\mathbb{Q}(u^{2}, t)$. In fact all the entries of the matrix
representatives for $G$ lie in the field $\mathbb{Q}(u^{2}, t)$.
This real field will be called the {\em invariant field
of definition} of $\Delta(u^{2},2t)$ as it is the most germane field
for our considerations. In particular, the cusp set of
$\Delta(u^{2},2t)$ can clearly be no larger than the field
$\mathbb{Q}(u^{2}, t) \cup \infty$

The main result of \cite{LR} is that there are rational choices of
parameters $(u^{2},2t)$ which give rise to nonarithmetic groups whose
cusp sets are precisely the rationals. Such groups we call {\em
pseudomodular}. There is a good deal of evidence that such groups
exist for fields more general than the rationals, that is to say,
their cusp sets are equal to their invariant field of definition \--
such groups we will describe as {\em maximally cusped}. It is these
groups which we will use to construct Dedekind sums; since our family
includes the modular group, it will include a construction of the
classical Dedekind sum. In this note we will show
\begin{theorem}
\label{main}
Suppose that $\Delta$ as above has invariant field of definition 
 $K$ and is maximally cusped.
Then associated to $\Delta$ is a function
$$S_{\Delta}\co K\cup \infty \rightarrow K\cup \infty$$
\end{theorem}
Such functions we say are {\em generalized Dedekind sums}.


\section{The construction}
Following \cite{KM}, we first construct an analogue of the Rademacher
$\phi$--function. Fix one of the groups $\Delta(u^{2},2t)$ of
\cite{LR}; (at this stage it is not necessary that the group be
pseudomodular) and suppose that its invariant field of definition
is $K$.

All once-punctured tori are hyperelliptic so we can adjoin
to this group the orientation-preserving involution $\tau$ which
conjugates the generators to their inverses, to form a new discrete
group $\Gamma$. The surface $F = {\bf H}^2/\Gamma$ is a sphere with
three cone points of angle $\pi$ and a cusp. Note that as an element
of $\GL(2,\mathbb{R})$, $\tau$ is represented by
the matrix $\begin{pmatrix}0&2u\cr-2/u& 0\cr\end{pmatrix}$, so
$\tau(\infty) = 0$.

Following \cite{KM}, we define an area $2$--cocycle
$$ \epsilon\co \Gamma \times \Gamma \rightarrow \mathbb{Z}$$
by setting $\epsilon(A,B) = \area(\infty,A\infty,AB\infty)/\pi$ where
this area is to be regarded as oriented, it follows that $\epsilon$
takes on the values $0,\pm 1$.

Equivalently, one can usefully think of $\epsilon(A,B)$ as the sign of
$AB\infty - A\infty$, where this is to be interpreted as zero if
either term of the difference is infinite.

Notice that $\epsilon$ is a cocycle, because the coboundary
$$\delta\epsilon(A,B,C) = 
\epsilon(B,C)-\epsilon(AB,C)+\epsilon(A,BC)-\epsilon(A,B)$$
involves four triangular areas and the first has vertices $(\infty, 
B\infty, BC\infty)$ which has same oriented area as $(A\infty, 
AB\infty, ABC\infty)$, so that taken together with
other three this forms a tetrahedron, and hence the total area is $0$.

\begin{lemma}
\label{cochain}
There is a unique $K$--valued $1$--cochain $\Gamma \rightarrow K$
with coboundary $\epsilon$.
\end{lemma}
\begin{proof}
Note that $\Gamma \cong \mathbb{Z}/2 * \mathbb{Z}/2 * \mathbb{Z}/2$ so that 
$$H^{1}(\Gamma ; \mathbb{Z}) \cong 0$$
and
$$H^{2}(\Gamma ; \mathbb{Z}) \cong \mathbb{Z}/2 \oplus \mathbb{Z}/2 \oplus \mathbb{Z}/
2$$
since the integral homology of $\mathbb{Z}/2$ is zero in odd dimensions
and $\mathbb{Z}/2$ in even dimensions. For our purposes, we need only
use that $H^2(\Gamma; K) =H^1(\Gamma; K) =0$. The fact that
$H^2(\Gamma; K) = 0$ implies immediately the existence of a
$K$--valued $1$--cochain with coboundary $\epsilon$.

We prove uniqueness as follows. If $\delta(\phi_1) = \epsilon =
\delta(\phi_2)$, then $\phi_i$ are both cocycles and hence since
$H^1(\Gamma; K ) = 0$, both are coboundaries. It follows that
there is a $0$--cochain $\beta$ with $\delta(\beta) = \phi_1 - \phi_2$. 
We are computing group cohomology with trivial coefficients, so that
this coboundary map is zero and $\phi_1 = \phi_2$ as required. 
\end{proof}
{\bf Definition}\qua We shall denote this $K$--valued $1$--cochain by $\phi$.


\subsection{Computation of $\phi$}
It will be useful to have a computation of the cochain $\phi$.
A consequence of Lemma \ref{cochain} is that there is a function 
$\phi\co \Gamma \rightarrow K$
which satisfies 
\renewcommand{\theequation}{$*$}
\begin{equation}
\phi(AB)-\phi(A)-\phi(B) = -\lambda \sign(AB\infty - A\infty) 
\end{equation}
\renewcommand{\theequation}{\arabic{equation}}
for some $\lambda \in K$ which will be determined. 

Taking $A = B = I$ we see that $\phi(I)=0$. Taking $A = B = -I$, we
also get $\phi(-I) = 0$. Taking $A=-I$ and $B=g$, we deduce from
$(*)$ that $\phi(g) = \phi(-g)$ for every $g \in \Gamma$.

More generally, if $A$ and $B$ both stabilise $\infty$, then the
relation says
$$\phi(AB) - \phi(A) - \phi(B)= 0$$
that is to say, $\phi$ is a homomorphism on $\stab(\infty)$. 

Note that in the group $\Gamma$, we have that $g_{1}g_{2}^{-1}\tau$ 
stabilises infinity and one checks easily that this is the generating 
matrix for the parabolic subgroup and is given by
$ \left( \begin{array}{cc} 1& t \\ 0& 1 \end{array}\right)$.

By scaling by an appropriate element of $K$, we may assume that
$\lambda$ is chosen so that $\phi$ maps this generating parabolic
matrix to $t$, so that $\phi$ is now determined on the parabolic subgroup.

It also follows from $(*)$ that
$$\phi(\alpha^{-1}) = - \phi(\alpha) = - \phi(-\alpha)$$
for any element $\alpha$, in particular, if $\xi$ is any projective 
involution in $\Gamma$,
(that is to say $ \xi^{2} = \pm I$)
we deduce that $\phi(\xi)=0$.

Now in the notation introduced above we have
$$\phi(g_1\tau) - \phi(g_1) - \phi(\tau) = -\lambda \sign(u^2-t+1)$$
Since $g_1\tau$ and $\tau$ are both projective involutions and recalling 
that
the groups in question are required to have $0 > 1+u^2-t$ we get 
$$ \phi(g_1) = -\lambda $$
By considering $\tau g_2$, a similar computation also shows $\phi(g_2) = -
\lambda $.


Now for any $ k \in K$, for which the matrix $\left( \begin{array}{cc}
1 & k \\ 0 & 1 \end{array} \right)$ lies in $\Gamma$, we have that
$$\phi\left( \left( \begin{array}{cc} 1 & k \\ 0 & 1 \end{array}
  \right) \tau \right) -
  \phi\left( \left( \begin{array}{cc} 1 & k \\ 0 & 1\end{array} \right)
  \right) - 0 = 0$$
In the special case that $k=t$, the leftmost term is the product 
$(g_{1}g_{2}^{-1}\tau)\tau = -g_{1}g_{2}^{-1}$, so we deduce from
the properties described above that $\phi(g_{1}g_{2}^{-1}) = t$.

Since $\left( \begin{array}{cc} t/u & -u \\ 1/u & 0 \end{array} 
\right) = g_{1}g_{2}^{-1}$, (or from purely geometric 
considerations) we see that $g_{2}g_{1}^{-1}\infty = 0$. 
Finally, noting that $g_{2}\infty = u^{2} > 0$ together with the 
relation
$$\phi(g_{2}g_{1}^{-1}) -\phi(g_2)+\phi(g_1) =-\lambda \sign(0 - u^{2})$$
it follows that $\lambda = -t$, since the leftmost term is $-t$ by the
previous calculation and the inverse rule. 


To sum up, we now have a complete inductive description of $\phi$ on
the group $\Gamma$, namely it satisfies
$$\phi(AB)-\phi(A)-\phi(B) = t \sign(AB\infty - A\infty) $$
and 
$$ \phi(g_1) = \phi(g_2) = \phi(g_{2}g_{1}^{-1}) = t $$
{\bf Remark}\qua This is in keeping with the computations of \cite{KM} 
which are for the modular group and have $\lambda = -3$.

\subsection{Generalized Dedekind sums}
Now fix some maximally-cusped $\Delta = \Delta(u^{2},2t)$ defined over 
the field $K$.

For any $M \in \Delta$, by applying the cocycle condition we have 
$$\phi\left(M\left( \begin{array}{cc} 1 & k \\ 0 & 1 \end{array}
  \right)\right) - \phi(M) - k = t\cdot \sign( M\infty - M\infty) = 0$$
from which it follows that
\renewcommand{\theequation}{$**$}
\begin{equation}
\phi\left(M\left( \begin{array}{cc} 1 & k \\ 0 & 1 \end{array}
  \right)\right) = \phi(M) + k
\end{equation}
\renewcommand{\theequation}{\arabic{equation}}
For $M\in \Delta\setminus \stab(\infty)$, set
$$\chi(M) = (M_{1,1} + M_{2,2})/M_{2,1}.$$ Since $M_{2,1}\neq 0$, the
value $\chi(M)$ is an element of the field $K$, since the groups
$\Delta$ consist of matrices of the shape $\sqrt r X$ for a matrix $X
\in \GL(2,K)$ and $r \in K$. Now a matrix computation shows that
$$\chi(M.\left( \begin{array}{cc} 1 & k \\ 0 & 1 \end{array} \right)) = 
\chi(M) + k$$
so that by taking the difference between this and $(**)$ we get
$$\phi\left(M\left( \begin{array}{cc} 1 & k \\ 0 & 1 \end{array}
  \right)\right) - \chi\left(M\left( \begin{array}{cc} 1 & k \\ 0 & 1
  \end{array} \right)\right) = \phi(M) - \chi(M)$$
which is to say the function
$$ S(M) = \phi(M) - \chi(M)$$
is invariant under right muliplication by the parabolic subgroup.

These observations are independent of whether $\Delta$ is maximally-cusped
or not. If we now assume that it is, we can define a {\em generalized Dedekind 
sum} as follows. 

Given any element $\kappa \in K$, since $\Delta$ is maximally cusped,
there is an element $M \in \Delta$ with $M(\infty) = \kappa$ and we may
set
$$S_{\Delta}(\kappa) = S(M)$$
The ambiguity in such $M \in \Delta$ is accounted for by right
multiplication by elements of the parabolic subgroup $\stab(\infty)$ so
that this function depends only on $\kappa$. We will define
$S_\Delta(\infty) = \infty$, and this defines the advertised function
in Theorem \ref{main}.

{\bf Remark}\qua
This construction gives a scalar multiple of the classical Dedekind sum
when $(u^{2},2t) = (1,6)$ (see \cite[section 0.8]{KM}).

{\bf Examples}\qua
It is proved in \cite{LR} that the group $\Delta(3/5,4)$ is
pseudomodular, so provides an example of a generalized Dedekind sum of
this type. It is not difficult to write a computer program which
computes its values based upon the iterative procedure outlined above.
A table of the groups currently proven to be pseudomodular (and some
conjectural examples) is provided in \cite{LR}. 

In subsequent work,
the authors have extended this table of conjectural examples to groups
which are maximally cusped for real quadratic number fields, for
example $\Delta(1, 2((1+\sqrt{13})/2))$ appears to be maximally cusped.
Questions about whether there are analogues of, for example, Dedekind
reciprocity and formulae of the classical type seem interesting and
appear worthy of further investigation.

\rk{Acknowledgements}This work was partially supported by the NSF and the
second author was partially supported by a grant from
the Texas Advanced Research Program.

\np

\Addresses

\end{document}